\documentclass[12pt]{article}
\usepackage{latexsym, amsmath, amssymb, amsfonts, amscd}
\usepackage{color}
\usepackage{graphics}
\usepackage[dvips]{epsfig}

\newtheorem{theorem}{Theorem}
\newtheorem{proposition}{Proposition}
\newtheorem{lemma}{Lemma}

\newtheorem{definition}{Definition}

\newtheorem{question}{Question}

\newcommand{\adL}{\mbox{\rm ad}_{\Lambda}}

\newcommand{\ad}{\mbox{\rm ad}}

\def\lcf{\lbrack\! \lbrack}
\def\rcf{\rbrack\! \rbrack}

\def\dbar{\overline\partial}

\newcommand{\CC}{\mathbb{C}}

\def\dbar{\overline\partial}

\def\oomega{\overline\omega}

\newcommand{\lieg}{\mathfrak{g}}
\newcommand{\liet}{\mathfrak{t}}
\newcommand{\liec}{\mathfrak{c}}
\newcommand{\lieh}{\mathfrak{h}}

\newcommand{\lbra}[2]{\lcf #1, #2 \rcf}

\newcommand{\bproof}{\noindent{\it Proof: }}
\newcommand{\eproof}{\hfill \qed \vspace{0.2in}}
\def\qed{\rule{2.3mm}{2.3mm}}

\begin{document}
\title{\bf Degeneracy of Holomorphic Poisson Spectral Sequence}
\author{
Yat Sun  Poon\thanks{ Address:
    Department of Mathematics, University of California at Riverside,
    Riverside, CA 92521, U.S.A.. E-mail: ypoon@ucr.edu.} }
%\date{August 24, 2016}
\maketitle
\begin{abstract} Through the theory of Lie bi-algebroids and generalized complex structures,
one could define a cohomology theory naturally associated to a holomorphic Poisson structure. It is known that
it is the hypercohomology of a bi-complex such that one of the two operators is the classical
$\dbar$-operator. Another operator is the adjoint action of the Poisson bivector with respect to
the Schouten-Nijenhuis bracket. The hypercohomology
is naturally computed by one of the two associated spectral sequences.
In a prior publication, the author of this article and his collaborators investigated
the degeneracy of this spectral sequence on the second page. In this note, the author investigates the
conditions for which this spectral sequence degenerates on the first page. Particular effort is devoted
to nilmanifolds with abelian complex structures.
\end{abstract}

%\tableofcontents

\section{Introduction}

It is well known that complex structures and symplectic structures are examples of
generalized complex structures in the sense of Hitchin
{\cite{Marco, Hitchin-Generalized CY, Hitchin-Instanton}}.
It is now also known that
holomorphic Poisson structure plays a fundamental role in generalized geometry
\cite{Bailey}. Since a key feature of generalized geometry is to put both complex
structures and symplectic structures within a single framework,
its deformation theory is of great interests
\cite{Goto1} \cite{Goto2} \cite{Hitchin-holomorphic Poisson}.
An understanding of
 cohomology theory of generalized geometry in general and
 holomorphic Poisson structures in particular becomes necessary.
On algebraic surfaces, there has been work done by Xu and collaborators
\cite{Hong-Xu} \cite{Xu}. For nilmanifolds with abelian
complex structures, there are recent work done by this author and his collaborators \cite{CFP} \cite{CGP} \cite{GPR}.
A common feature of these work is to recognize the cohomology as the hypercohomology of a bi-complex.

In computation of hypercohomology of bi-complex,  theoretically one could apply
one of the two naturally defined spectral sequences to
complete the task. In the case of a holomorphic Poisson structure, the first page of
one of the spectral sequences consists of the Dolbeault cohomology of a complex manifold with coefficients in the
sheaf of germs of holomorphic polyvector fields. As the Dolbeault cohomology is a well known classical object,
it is natural to determine when this spectral
sequence degenerates fast. In \cite{CFP} and \cite{CGP}, we have seen situation when the spectral sequence
degenerates on the second page. We push our analysis in the past and investigate the possibility when
degeneracy on the first page occurs.

In particular, as a corollary of our Theorem \ref{main theorem} in this article,
we could formulate a result as below.

\begin{theorem}\label{corollary theorem}
Let $M=G/\Gamma$ be a 2-step nilmanifold with abelian complex structure.
Let $\liec$ be the center of the Lie algebra $\lieg$ of the simply connected Lie group $G$, and
$\liet$ its quotient.
Let $\lieg^{1,0}=\liet^{1,0}\oplus \liec^{1,0}$ be the space of invariant $(1,0)$-vectors.
Assume that $\dim\liec^{1,0}=1$. Let $\overline\rho$ span $\liec^{*(0,1)}$.
Suppose that $\Lambda_1$ is in $\liet^{1,0}\otimes \liec^{1,0}$.
Then,
\begin{itemize}
\item $\Lambda_1$ is a holomorphic Poisson structure.
 \item The spectral sequence of its bi-complex degenerates on the second page.
\item If in addition, $d\rho$ is a non-degenerate $(1,1)$-form, the spectral sequence
degenerates on the first page.
\end{itemize}
\end{theorem}

We will develop a proof leading to the above result and provide examples to illustrate whether one could
relax the constrains to attain the same result.

\section{Complex and Generalized Complex Structures}

In this section, we review the basic background materials as seen in \cite{CGP} to set up notations.

Let $M$ be a smooth manifold. Denote its tangent bundle by $TM$ and the co-tangent bundle by $T^*M$.
A generalized complex structure on an even-dimensional manifold $M$  \cite{Marco} \cite{Hitchin-Generalized CY}
 is a
subbundle $L$ of the direct sum
${\cal T}=(TM\oplus T^*M)_{\CC}$ such that
 \begin{itemize}
 \item  $L$ and its conjugate bundle $\overline{L}$ are transversal;
 \item  $L$ is maximally isotropic with respect to the natural pairing on ${\cal T}$;
  \item and the space of sections of $L$ is closed with respective to the Courant bracket.
  \end{itemize}

Given a generalized complex structure,  the pair of bundles $L$
and $\overline L$ makes a (complex) Lie bi-algebroid. The composition of the
 inclusion of $L$ and $\overline L$ in
${\cal T}$ with the natural projection onto the summand $TM_{\CC}$ becomes the anchor map
of these Lie algebroids. It is denoted by $\phi$.
Via the canonical non-degenerate pairing on the bundle $\cal T$,
 the bundle $\overline L$ is complex linearly identified to the dual of $L$.
 Therefore, the Lie algebroid differential of $\overline L$ acts on $L$. It extends to a differential on
 the exterior algebra of $L$.
 For the calculus of Lie bi-algebroids,
 we follow the conventions in \cite{Mac}. In particular,
 for any element $\Gamma$ in $C^\infty(M, \wedge^kL)$ and elements ${a}_1, \dots,
 {a}_{k+1}$ in $C^\infty(M, {\overline{L}})$, the Lie algebroid differential  of $\Gamma$  is defined by
 the Cartan Formula as in exterior differential algebra, namely
 \begin{eqnarray}
&& (\dbar_{L}\Gamma)({a}_1, \dots, {a}_{k+1})
=\sum_{r=1}^{k+1}(-1)^{r+1}\phi(a_r)(\Gamma({a}_1, \dots,{\hat{a}}_{r} ,\dots, {a}_{k+1}))\nonumber\\
&\quad& \quad +
\sum_{r<s}(-1)^{r+s}\Gamma(\lbra{a_r}{a_s}, {a}_1, \dots,{\hat{a}}_{r} ,\dots,{\hat{a}}_{s} ,\dots,{a}_{k+1}).
\label{algebroid differential}
\end{eqnarray}
The space of sections of the bundle $\overline{L}$ is closed if and only if $\dbar_{L}\circ \dbar_{L}=0$.

Typical examples of generalized complex structures
 are classical complex structures and symplectic structures on a manifold.
 In this section, we focus on the former.

Let $J: TM\to TM$ be an integrable complex structure on the manifold $M$.
The complexified tangent bundle $TM_\CC$ splits into the direct sum of
bundle of $(1,0)$-vectors $TM^{1,0}$ and bundle of $(0,1)$-vectors $TM^{0,1}$.
Their $p$-th exterior products are respectively denoted by $TM^{p,0}$ and $TM^{0,p}$.
Denote their dual bundles by $TM^{*(p,0)}$ and $TM^{*(0,p)}$ respectively.

Define $L= TM^{1,0}\oplus TM^{*(0,1)} $. One gets a generalized complex structure.
Its dual is its complex conjugate ${\overline{L}}=TM^{0,1}\oplus TM^{*(1,0)}$.
When one restricts the Courant bracket from the ambient bundle
${\cal T}=(TM\oplus T^*M)_{\CC}$ to the subbundles $L$ and $\overline{L}$, then one recovers the
Schouten-Nijenhuis bracket or simply known as the Schouten bracket in classical deformation. The
Schouten bracket between
$(1,0)$-vector fields is the Lie bracket of vector fields;
the Schouten bracket between a $(1,0)$-vector field and a $(0,1)$-form is the Lie derivative of a form by a
vector field. The Schouten bracket between two $(0,1)$-forms is equal to zero. These brackets are
extended to higher exterior product by observing the rule of exterior multiplication \cite{Mac}.

With respect to the Lie algebroid $\overline{L}$,
we get its differential $\dbar$  as defined in
(\ref{algebroid differential}).
\begin{equation}
\dbar: C^\infty(M, L ) \to
C^\infty(M, \wedge^2L).
\end{equation}
It is extended to a differential of exterior algebras:
\begin{equation}
\dbar: C^\infty(M, \wedge^pL) \to
C^\infty(M, \wedge^{p+1}L).
\end{equation}
It is an elementary exercise in computation of Lie algebroid differential that when $\dbar$ is restricted to
$(0,1)$-forms, it is the classical $\dbar$-operation in complex manifold theory; and
\[
\dbar: C^\infty(M, TM^{*(0,1)}) \to C^\infty(M, TM^{*(0,2)})
\]
is the $(0,2)$-component of the exterior differential \cite{Brian}.
Similarly, when the Lie algebroid differential is restricted to
$(1,0)$-vector fields, then
\[
\dbar:  C^\infty(M, TM^{1,0} )\to C^\infty(M, TM^{*(0,1)}\otimes TM^{1,0});
\]
and it is the Cauchy-Riemann operator as seen in \cite{Gau}.

By virtue of $L$ and $\overline{L}$ being a pair of Lie bi-algebroid, the
space $C^\infty (M, \wedge^\bullet L)$ together with the Schouten bracket,
exterior product and the Lie algebroid differential $\dbar$ form a differential Gerstenhaber
algebra  {\cite{Mac}} \cite{Poon}. In particular, if $a$ is a smooth section
of  $\wedge^{|a|}L$
and $b$ is a smooth section of $\wedge^{|b|}L$, then
\begin{eqnarray}
\dbar \lbra{a}{b} &=& \lbra{\dbar a}{ b}+(-1)^{|a|+1}\lbra{a}{\dbar b};\\
\dbar (a\wedge b) &=& (\dbar a)\wedge b+(-1)^{|a|}a\wedge (\dbar b),
\end{eqnarray}

Since $\dbar\circ\dbar=0$, one obtains the Dolbeault cohomology with coefficients in holomorphic
polyvector fields. Denote the sheaf of germs of sections of
the $p$-th exterior power of the holomorphic tangent bundle by $\Theta^p$, we have
\[
H^\bullet(M, \wedge^\bullet \Theta)\cong
\oplus_{p,q\geq 0}H^q(M, \Theta^p).
\]
In subsequent computation, when $p=0$, $\Theta^p$ means to represent the structure sheaf $\mathcal{O}$ of the
complex manifold $M$.

Due to the compatibility between $\dbar$ and
 the Schouten bracket $\lbra{-}{-}$ and the compatibility between $\dbar$ and the exterior product $\wedge$
 as noted above, the Schouten bracket and exterior product descend to
the cohomology space $H^\bullet(M, \Theta^\bullet)$. In other words,
the triple $(H^\bullet(M, \wedge^\bullet TM^{1,0}), \lbra{-}{-}, \wedge)$ forms a Gerstenhaber algebra.
 When we ignore the exterior product,  we call it  a Schouten algebra.
 For example, by center of the Schouten algebra, we mean the collection of  elements $A$ in
 $H^\bullet(M, \Theta^\bullet)$
  such that $\lbra{A}{B}=0$ for all $B$ in this space.

\section{Holomorphic Poisson Bi-Complex}

A holomorphic Poisson structure on a complex manifold $(M, J)$ is a holomorphic bi-vector field $\Lambda$
such that  $\lbra{\Lambda}{\Lambda}=0$.
The corresponding bundles as generalized complex structure for ${\Lambda}$ are the pair of
bundles of graphes $L_{\overline{\Lambda}}$ and ${\overline L}_{\Lambda}$ where
\begin{equation}
{\overline L}_{\Lambda}=\{{\overline\ell}+\Lambda(\overline{\ell}): {\overline\ell}\in {\overline L}\}.
\end{equation}

 While the pair of bundles
 $L_{\overline{\Lambda}}$ and ${\overline L}_{\Lambda}$ form naturally a Lie bi-algebroid,
 so does the pair $L$ and ${\overline L}_{\Lambda}$ \cite{LWX}. From this perspective, the Lie
 algebroid differential of the deformed generalized complex structure ${\overline L}_{\Lambda}$ acts on
 the space of sections of the bundle $L$.

Any smooth section of the bundle $L$ is the sum of a section $v$ of $T^{1,0}M$ and a section
$\oomega$ of  $T^{*(0,1)}M$.
Given a holomorphic Poisson structure $\Lambda$, define $\adL$ by
\begin{equation}
\adL(v+\oomega)=\lbra{\Lambda}{v+\oomega}.
\end{equation}

 \begin{proposition}\label{d lambda}\mbox{\rm See (Proposition 1) in \cite{GPR}.} The action
 of the Lie algebroid differential
 of ${\overline L}_{\Lambda}$  on $L$ is given by
 \begin{equation}
 \dbar_{\Lambda}=\dbar+\adL: C^\infty(M, L)\to
 C^\infty(M, \wedge^2L).
 \end{equation}
 \end{proposition}

 The operator $\dbar_{\Lambda}$ extends to act on the exterior algebra of $T^{1,0}M\oplus T^{*(0,1)}M=L$.
  From now on, for $n\geq 0$ denote
 \begin{equation}
 K^n=C^\infty (M, \wedge^n L).
 \end{equation}
 For $n< 0$, set $K^n=\{ 0\}$.

 Since $\Lambda$ is a holomorphic Poisson structure, the closure of the space of section of the
 corresponding bundles ${\overline L}_{\Lambda}$ is equivalent to  $\dbar_{\Lambda}\circ \dbar_{\Lambda}=0$.
 Therefore, one has a complex with $\dbar_{\Lambda}$ being a differential.
 \begin{definition}
For all $n\geq 0$, the $n$-th holomorphic Poisson cohomology of the holomorphic Poisson structure $\Lambda$ is the space
 \begin{equation}
 H^n_{\Lambda}(M):=\frac{{\mbox{\rm kernel of  }}\ \dbar_\Lambda: K^n \to K^{n+1} }
 {{\mbox{\rm image of }}\ \dbar_\Lambda: K^{n-1}\to K^n}.
 \end{equation}
 \end{definition}
Given Proposition \ref{d lambda}, the identity
 $\dbar_{\Lambda}\circ \dbar_{\Lambda}=0$ is equivalent to a system of three.
 \begin{equation}
 \dbar\circ\dbar=0, \quad
 \dbar\circ \adL+\adL\circ \dbar=0,
 \quad
 \adL\circ\adL=0.
 \end{equation}
 The first identity is equivalent to the complex structure $J$ being integrable;
 the second identity is equivalent to $\Lambda$ being holomorphic, and
 the third is equivalent to $\Lambda$ being Poisson.
Define $A^{p,q}=C^\infty(M, TM^{p,0}\otimes TM^{*(0,q)})$, then
\begin{equation}
\adL: A^{p,q}\to A^{p+1, q},
\quad
\dbar: A^{p,q} \to A^{p, q+1};
\quad
\mbox{ and }
\quad
K^n=\oplus_{p+q=n}A^{p,q}.
\end{equation}
Therefore, we obtain a bi-complex. We arrange the double indices $(p,q)$
in such a way that $p$ increases horizontally so that
$\adL$ maps from left to right, and
$q$ increases vertically so that $\dbar$ maps from bottom to top.
\begin{definition} Given a holomorphic Poisson structure $\Lambda$, its Poisson bi-complex is
the triple $\{ A^{p,q}, \adL, \dbar\}$.
\end{definition}
It is now obvious that the (holomorphic) Poisson cohomology $H^\bullet_\Lambda(M)$
theoretically could be computed by each one of the
two naturally defined spectral sequences.
We choose a filtration given by
\[
F^pK^n=\oplus_{p'+q=n, p'\geq p}A^{p',q}.
\]
The
lowest differential is $\dbar: A^{p,q} \to A^{p, q+1}$. Therefore, the first sheet of the
spectral sequence is the Dolbeault cohomology
\begin{equation}
E_1^{p,q}=H^q(M, \Theta^p).
\end{equation}
The first page
of the spectral sequence is given as below.
\[
\begin{array}
[c]{cccccccc}
& \cdots & \rightarrow & \cdots & \rightarrow & \cdots & \rightarrow & \\
& H^{n}(M,\mathcal{O}) & \rightarrow & H^{n}(M,\Theta) & \rightarrow
& H^{n}(M,\Theta^{2}) & \rightarrow & \\
%&  &  &  &  &  &  & \\
& \cdots & \rightarrow & \cdots & \rightarrow & \cdots & \rightarrow & \\
& H^{1}(M,\mathcal{O}) & \rightarrow & H^{1}(M,\Theta) & \rightarrow
& H^{1}(M,\Theta^{2}) & \rightarrow & \\
%&  &  &  &  &  &  & \\
& H^{0}(M,\mathcal{O}) & \rightarrow & H^{0}(M,\Theta) & \rightarrow
& H^{0}(M,\Theta^2) & \rightarrow &
%\\ E_{1}^{p,q} &  & {\mbox{\rm{ad}}}_{\Lambda}=[\Lambda,-] &  &  &  &  &
\end{array}
\]
The differential on this page is
 \begin{equation}\label{d1}
 d_1^{p,q}=\adL: H^q(M, \Theta^p)\to H^q(M, \Theta^{p+1}).
 \end{equation}

\begin{question} When will the spectral sequence of a Poisson bi-complex
degenerates on the first page? In other words, when will
$\adL \equiv 0$ for all $p,q$?
\end{question}

Prior to us investigating non-trivial sufficient conditions for spectral sequence to
degenerate on the first page, we consider necessary conditions.
For instance, we should derive the necessary conditions for
the $d_1$-map on the first row and the first column vanishes.

Given the definition of $d_1$ in (\ref{d1}), the following observation
regarding the first row of the first page in the spectral sequence is trivial.

\begin{proposition}\label{polyvector fields} Suppose that the bi-complex
of a holomorphic Poisson structure $\Lambda$ degenerates on the first page,
then for any holomorphic polyvector fields $\Upsilon$, i.e. elements in  $\oplus_{p\geq 0}H^0(M, \Theta^p)$,
$\lbra{\Lambda}{\Upsilon}=0$.
\end{proposition}

Note that the restriction of the Schouten bracket from the full cohomology space $\oplus_{p,q}H^q(M, \Theta^p)$
to the subspace of holomorphic polyvector fields $\oplus_{p\geq 0}H^0(M, \Theta^p)$ turns the latter into
a subalgebra. The above proposition means that  $\Lambda$ is in the
 center  of the Schouten algebra $\oplus_{p\geq 0}H^0(M, \Theta^p)$.

The observation above immediately demonstrates that the spectral sequence of
well known holomorphic Poisson structures will not degenerate on its first page.
Complex projective spaces are the cases at hand. On the other hand, if the space of holomorphic
polyvector fields is abelian with respect to the Schouten bracket,
then the obstruction from the first row to degenerate vanishes.

Next, we consider the $d_1$-map on the first column.
\begin{equation}
 d_1^{0,q}=\adL: H^q(M, {\cal O})\to H^q(M, \Theta).
 \end{equation}
Suppose $\oomega$ represents a class in $H^q(M, \cal{O})$, then $\adL(\oomega)$ represents the zero class in
$H^q(M, \Theta)$ if it is $\dbar$-exact. Therefore, we have the following observation.

\begin{proposition} Suppose that the bi-complex of a holomorphic
Poisson structure $\Lambda$ degenerates on the first page,
then for any $\dbar$-closed $(0,q)$-form $\oomega$, there exists $\Gamma\in
C^\infty(M, T^{1,0}M\otimes T^{*(0,q-1)})$ such that
\begin{equation}
\adL(\oomega)=\dbar(\Gamma).
\end{equation}
\end{proposition}

Given the obvious necessary conditions by observing the $d_1$-map on the first row and the
first column on the first page, we now could post a question to guide our investigation.
\begin{question} Suppose that $\Lambda$ is in  the center of the Schouten algebra
of the holomorphic polyvector fields $\oplus_{p\geq 0}H^0(M, \Theta^p)$ and that
for any $\dbar$-closed $(0,q)$-form $\oomega$,
$\adL(\oomega)$ is $\dbar$-exact. Is it necessarily that $d_1^{p,q}\equiv 0$ for all $p,q$?
\end{question}

In the rest of this paper, we find answers to these questions on a class of nilmanifolds.

\section{2-step Nilmanifolds}

A compact manifold $M$ is a nilmanifold if there exists a
simply-connected nilpotent Lie group $G$
and a lattice subgroup $\Gamma$ such that $M$ is diffeomorphic to $G/\Gamma$.
We denote the Lie algebra of the group
$G$ by $\lieg$ and its center by $\liec$. The step of the
nilmanifold is the nilpotence of the Lie algebra $\lieg$. A left-invariant complex
structure $J$ on $G$ is said to be abelian if on the Lie algebra $\lieg$, it satisfies the conditions
$J\circ J=-$identity and $\lbra{JA}{JB}=\lbra{A}{B}$ for all $A$ and $B$ in the Lie algebra $\lieg$.
If one complexifies the algebra $\lieg$ and denotes the $+i$ and $-i$ eigen-spaces of $J$ respectively
by $\lieg^{1,0}$ and $\lieg^{0,1}$, then the invariant complex structure $J$ is abelian
if and only if the complex Lie algebra $\lieg^{1,0}$ is abelian.

Denote $\wedge^k\lieg^{1,0}$ and $\wedge^k{\lieg}^{*(0,1)}$ respectively by
$\lieg^{k,0}$ and ${\lieg}^{*(0,k)}$. We will use the following notation.
\[
B^{p,q}=\lieg^{p,0}\otimes {\lieg}^{*(0,q)}.
\]

Assume that the Lie algebra $\lieg$ is 2-step nilpotent, i.e. $[\lieg,\lieg]\subset \liec$.
In such a case, we call the manifold $M=G/\Gamma$ a 2-step nilmanifold \cite{MPPS-2-step}.

On the nilmanifold $M$, we consider $\lieg^{k,0}$ as invariant $(k,0)$-vector fields and
${\lieg}^{*(0,k)}$ as invariant $(0,k)$-forms. It yields an inclusion map
\[
B^{p,q}
\hookrightarrow A^{p,q}=C^\infty(M, T^{p,0}M\otimes T^{*(0,q)}M).
\]
When the complex structure is also invariant, $\dbar$ sends
$B^{p,q}$ to
$B^{p,q+1}$.
Given an invariant complex structure and an invariant holomorphic Poisson structure $\Lambda$,
$\adL$ sends
$B^{p,q}$ to
$B^{p+1,q}$.
Restricting $\dbar$ to $B^{p,q}$, we then consider the invariant cohomology.
\begin{equation}\label{hqp}
H^q( \lieg^{p,0})=
\frac
{\mbox{kernel of } \dbar: B^{p,q}\to B^{p,q+1}}
{\mbox{image of } \dbar: B^{p,q-1}\to B^{p,q}}.
\end{equation}
The inclusion map yields a homomorphism of cohomology:
\[
H^q({\lieg}^{p,0})\hookrightarrow H^q(M, \Theta^p).
\]

\begin{theorem}{\rm (See \cite{CFP} and \cite{CGP})}\label{quasi isomorphic}
On a 2-step nilmanifold $M$ with an invariant
abelian complex structure, the inclusion $B^{p,q}$ in
$A^{p,q}=C^\infty(M, T^{p,0}M\otimes T^{*(0,q)}M)$ induces an isomorphism of cohomology. In other words,
\[
H^q({\lieg}^{p,0})\cong H^q(M, \Theta^p).
\]
\end{theorem}

Given any element $A$ in $B^{p,q}$ , it acts on
$B^{p,q}$ by the Schouten bracket. We denote its action by $\ad_{A}$. i.e.
\[
\ad_A(B)=\lbra{A}{B}.
\]
An element $A$ is in the center  of the Schouten algebra
$\oplus_{p,q}B^{p,q}$ with respect to the Schouten bracket
$\lbra{-}{-}$ if and only if
$\ad_A\equiv 0$.
Similarly, an element $A$ in $H^q({\lieg}^{p,0})$ is in the center of the Schouten algebra
$\oplus_{p,q}H^q({\lieg}^{p,0})$ if $\ad_A(B)$ is equal to zero on the cohomology level for
any $B$ in $\oplus_{p,q}H^q({\lieg}^{p,0})$.

\

 Let  $\liet=\lieg/\liec$.
 Below are some facts shown in Sections 2 and 3 of \cite{MPPS-2-step}.
 Since $\lieg$ is 2-step nilpotent, $\liet$ is abelian.
 As a vector space, $\lieg^{1,0}=\liet^{1,0}\oplus \liec^{1,0}$,
 and  $\lieg^{*(1,0)}=\liet^{*(1,0)}\oplus \liec^{*(1,0)}$.
 The only non-trivial Lie brackets in $\lieg^{1,0}\oplus \lieg^{0,1}$
 are of the form
 $ [\liet^{1,0},\liet^{0,1}]\subset \liec^{1,0}\oplus \liec^{0,1}.$

Explicitly, there exists a real basis $\{X_k,JX_k:1\leq k\leq n\}$ for $\liet$ and
$\{Z_\ell, JZ_\ell:1\leq \ell \leq m\}$ a real basis for $\liec$. The corresponding complex bases
for $\liet^{1,0}$ and $\liec^{1,0}$ are respectively composed of the following elements:
\begin{equation}
T_k=\frac{1}{2}(X_k-iJX_k) \quad \mbox{ and } \quad W_\ell=\frac{1}{2}(Z_\ell-iJZ_\ell).
\end{equation}
The structure equations of $\lieg$ are determined by
\begin{equation}\label{structure eq}
\lbra{\overline{T}_k}{T_j}=\sum_{\ell}E_{kj}^\ell W_\ell-\sum_\ell{\overline{E}}_{jk}^\ell{\overline{W}}_\ell
\end{equation}
for some constants $E_{kj}^\ell$.
Let $\{\omega^k: 1\leq k\leq n\}$ be the dual basis for $\liet^{*(1,0)}$, and let
$\{\rho^\ell: 1\leq \ell\leq m\}$ be the dual basis for $\liec^{*(1,0)}$.
The dual structure equations for (\ref{structure eq}) are
\begin{equation}\label{dual 1}
d\rho^\ell=\sum_{i,j}E_{ji}^\ell\omega^i\wedge\oomega^j \quad \mbox{ and }
\quad d\omega=0.
\end{equation}
Equivalently,
\begin{equation}\label{dual 2}
d\overline{\rho}^\ell=-\sum_{i,j}\overline{E}^\ell_{ji}\omega^j\wedge\oomega^i
\quad \mbox{ and }
\quad d\oomega=0.
\end{equation}
It follows that
\begin{equation}\label{adj-T-on-rho-bar}
\lbra{T_j}{\overline{\rho}^\ell}={\cal L}_{T_j}{\overline{\rho}^\ell}
=\iota_{T_j}d{\overline{\rho}^\ell}=-\sum_i\overline{E}^\ell_{ji}\oomega^i.
\end{equation}

By Cartan Formula (\ref{algebroid differential}),
\begin{equation}
\dbar {T_j}=\sum_{k, \ell}E_{kj}^\ell\oomega^k\wedge W_\ell,
\end{equation}

The consequence of the above computation is twofold. One
is about the structure of Lie algebroid differential $\dbar$ on
$\wedge^\bullet(\lieg_\CC\oplus \lieg_\CC^*)$.
Another is about structure of Schouten bracket on the same space.
On differential, the above computation could be summarized as below.

\begin{lemma}\label{image of dbar}
$
\lieg^{*(0,1)}\oplus \liec^{1,0}\subseteq \ker\dbar$ and
$\dbar \liet^{1,0}\subseteq \liet^{*(0,1)}\otimes \liec^{1,0}$.
Moreover,
\[
\dbar (\liet^{k,0}\otimes \liec^{\ell,0}\otimes \liet^{*(0,a)}\otimes \liec^{*(0,b)})
\subseteq \liet^{k-1,0}\otimes \liec^{\ell+1,0}\otimes \liet^{*(0,a+1)}\otimes \liec^{*(0,b)}.
\]
\end{lemma}

In subsequent presentation, we suppress the notations for the vector spaces,
and simply keep track of the  quadruple of indices to indicate the composition of the components
involved.
In particular, every element in the Schouten algebra
$\oplus_{p,q}B^{p,q}$ decomposes into a sum of
different types according to their indices $(k,\ell; a,b)$, with $k+\ell=p$ and $a+b=q$.

With this notation, Lemma \ref{image of dbar} above is summarized as
\begin{equation}\label{dbar action}
\dbar(k, \ell; a, b)\subseteq (k-1, \ell+1; a+1, b).
\end{equation}
Note that in our notations, whenever anyone of the four indices $(k,\ell; a,b )$ is less than zero,
we mean the trivial vector space $\{0\}$.
Due to the assumption that $\dim\liec^{1,0}=1$, the space
$(k-1, \ell+1; a+1, b)$ is trivial when $\ell\geq 1$.

Next, we turn our attention to the structure of the Schouten bracket. In particular, we are
interested in the adjoint action of elements in $\liec^{1,0}$ and $\liet^{1,0}$.
Since the complex algebra $\lieg^{1,0}$ is abelian, it is clear that for any $W$ in $\liec^{1,0}$, the action
of $\ad_{W}$  on $\lieg^{1,0}$ is identically zero.
Since for any $\oomega \in \liet^{*(0,1)}$, $d\oomega=0$. It follows that $\ad_{W}\oomega=0$.
For any $\overline{\rho}\in \liec^{*(0,1)}$, $d\overline{\rho}\in \liet^{*(1,1)}$. it follows that
$\ad_{W}\overline{\rho}=0$ as well. Therefore, the action of
$\ad_{W}$ on $\lieg^{*(0,1)}$ is  identically zero; and subsequently $\ad_{W}\equiv 0$ on
$\lieg^{1,0}\oplus\lieg^{*(0,1)}$. For future reference we
summarize our observation on the central vector field $W$ in a lemma.

\begin{lemma}\label{central vector field} The space $\liec^{1,0}$ is in the center of
the Schouten algebra $(\oplus_{p,q\geq 0}B^{p,q}, \lbra{-}{-})$.
\end{lemma}

Similarly, for any $T\in \liet^{1,0}$, the restriction of $\ad_{T}$ on $\lieg^{1,0}\oplus \liet^{*(0,1)}$ is
identically zero, and the image of $\ad_{T}$ acting on $\liec^{*(0,1)}$ is contained in $\liet^{*(0,1)}$.

\begin{lemma}\label{adjoint T} For any $T\in \liet^{1,0}$,
$\lieg^{1,0}\oplus \liet^{*(0,1)}\subseteq \ker \ad_{T}$,  and
$\ad_T (\liec^{*(0,1)}) \subseteq \liet^{*(0,1)}.$
\end{lemma}

\begin{lemma}{\rm See also \cite[Lemma 7]{CGP}.} \label{image of adjoint}
 Suppose that $\Lambda_1=W\wedge T$ and $\Lambda_2=T_1\wedge T_2$ where
$W\in \liec^{1,0}$ and $T$, $T_1$, $T_2$ $\in \liet^{1,0}$, then
\begin{eqnarray}
\ad_{\Lambda_1} (k, \ell; a,b) &\subseteq & (k, \ell+1; a+1, b-1).
\label{image of adjoint 1}\\
\ad_{\Lambda_2} (k, \ell; a,b) &\subseteq & (k+1, \ell; a+1, b-1).
\label{image of adjoint 2}
\end{eqnarray}
\end{lemma}
\bproof
Whenever $V_1, V_2$ are in $\lieg^{1,0}$ and $\Phi$ in
$B^{p,q}$, then
\[
\lbra{V_1\wedge V_2}{\Phi}=V_1\wedge\lbra{ V_2}{\Phi}-V_2\wedge\lbra{ V_1}{\Phi},
\]
The lemma in question now follows the previous two.
\eproof

In particular, for any $\overline\rho\in \liec^{*(0,1)}$,
\begin{equation}
\lbra{\Lambda_1}{\overline{\rho}}\in
\liec^{1,0}\otimes \liet^{*(0,1)}
\quad \mbox{ and }
\quad
\lbra{\Lambda_2}{\overline{\rho}}\in  \liet^{1,0}\otimes\liet^{*(0,1)}.
\end{equation}

\section{Computation of the map $d_2$}

Now on a 2-step nilmanifold $M$ with abelian complex structure, if $\dim\liec^{1,0}=1$, then
\[
\lieg^{2,0}=(\liet^{1,0}\otimes\liec^{1,0})\oplus \liet^{2,0}.
\]
Therefore, an invariant bi-vector $\Lambda$ decomposes into a sum of two types.
\[
\Lambda=\Lambda_1+\Lambda_2
\]
where $\Lambda_1\in \liet^{1,0}\otimes\liec^{1,0}$ and
$\Lambda_2\in \liet^{2,0}$. Since the complex structure is abelian,
$\lbra{\Lambda}{\Lambda}=0$. Therefore, all bi-vectors are Poisson.
In subsequence computation, we assume that $\Lambda$ is  holomorphic.

For completeness we quickly review a computation of the map $d_2$  \cite{CGP}.
Recall that
\begin{equation}
E_2^{p,q}=
\frac
{\mbox{kernel of } \adL: H^q( \lieg^{p,0})\to H^q( \lieg^{p+1,0})}
{\mbox{image of } \adL: H^q( \lieg^{p-1,0})\to H^q( \lieg^{p,0})}.
\end{equation}

An element in $E_2^{p,q}$ is represented by an element  $\Upsilon\in \lieg^{p,0}\otimes\lieg^{*(0,q)}$ such that
it is in $H^q( \lieg^{p,0})$ and in the kernel of $\adL$. Equivalently,
$\Upsilon$ is $\dbar$-closed and $\adL\Upsilon$ is $\dbar$-exact. Therefore, there exists
$\Gamma \in \lieg^{p+1,0}\otimes\lieg^{*(0,q-1)}$ such that
\begin{equation}\label{potential}
\dbar\Upsilon =0 \quad {\mbox{ and }} \quad \adL\Upsilon=\dbar\Gamma.
\end{equation}
By definition,
$d_2[\Upsilon]$ is represented by $\adL\Gamma.$

When $\dim_{\mathbb{C}}\liec^{1,0}=1$,
the space $(k,\ell; a,b)$ is equal to the zero set except when $b\in\{0,1\}$ and $\ell\in\{0,1\}$.
Under this assumption, we consider the components of  $\Upsilon$.
By Lemma \ref{image of adjoint}, its component in $(k, \ell; a,0)$ is mapped to zero by $\adL$.
In such case, $\Gamma$ could be chosen to be zero,
and hence $d_2$ maps this component of $\Upsilon$ to  zero.

The non-trivial case is when $\Upsilon$ has a component in $(k,\ell; a,1)$.
By (\ref{image of adjoint 1})
and (\ref{image of adjoint 2}),
$\adL(\Upsilon)$ is contained in the direct sum of the following two types of subspaces
\[
(k, \ell+1; a+1, 0), \quad \quad (k+1, \ell; a+1, 0 ).
\]
By Lemma \ref{image of dbar} or equivalently Identity (\ref{dbar action}),
for $\Gamma$ to be a solution of (\ref{potential}), it must come from the
direct sum of the following two types of subspaces.
\[
 (k+1, \ell; a, 0), \quad \quad  (k+2, \ell-1; a, 0 ).
\]
By (\ref{image of adjoint 1}) and (\ref{image of adjoint 2}), elements of these types are in the kernel of
 $\ad_{\Lambda_1}$ and $\ad_{\Lambda_2}$. We conclude that $\adL\Gamma=0$, and hence $d_2[\Upsilon]=0$.
 So we recover a theorem in \cite{CGP}.

 \begin{theorem}{\rm \cite[Theorem 2]{CGP}}
 Suppose that $M$ is a 2-step nilmanifold with an abelian complex structure. Suppose that the
 center of the Lie algebra of the simply connected covering space is real two-dimensional, then for any
 invariant holomorphic Poisson structure, the spectral sequence of the Poisson bi-complex
  degenerates on the second page.
 \end{theorem}

\section{Computation of the map $d_1$}

Next, we push the work in the previous section to study first page degeneracy. First of all, we explore the necessary condition for
the spectral sequence of the bi-complex of
$\Lambda=\Lambda_1+\Lambda_2$   to degenerate on the first page.

Consider the adjoint action of $\Lambda$ on the zero-th row on the first page of the spectral sequence.
\[
\adL: H^0(\lieg^{p,0})\rightarrow H^0(\lieg^{p+1,0}).
\]
Since the complex structure is abelian, this map is identically zero.
As $\Lambda=\Lambda_1+\Lambda_2$, by linearity of the Schouten bracket
$\adL=\ad_{\Lambda_1}+\ad_{\Lambda_2}$.

Since $\overline{\rho}$ is $\dbar$-closed, it represents an element in $H^1(\mathbb{C})$.
If $\adL$ is identically zero, then in particular,
$\ad_{\Lambda_1}\overline{\rho}+\ad_{\Lambda_2}\overline{\rho}$ is $\dbar$-exact.

 With respect to the type decomposition $(k, \ell; a,b)$,
$\overline\rho$ is type $(0,0; 0,1)$. By (\ref{image of adjoint 1}),  $\ad_{\Lambda_1}{\overline\rho}$ is type
$(0,1; 1,0)$. By (\ref{image of adjoint 2}),  $\ad_{\Lambda_2}{\overline\rho}$ is type
$(1,0;1,0)$. As $\ad_{\Lambda_1}{\overline\rho}$ and $\ad_{\Lambda_2}{\overline\rho}$ are of different types,
$\ad_{\Lambda_1}\overline{\rho}+\ad_{\Lambda_2}\overline{\rho}$ is $\dbar$-exact if and only if
$\ad_{\Lambda_1}\overline{\rho}$ and $\ad_{\Lambda_2}\overline{\rho}$ are both $\dbar$-exact.

However, by Lemma \ref{image of dbar} or equivalently Identity (\ref{dbar action})
the image of the $\dbar$-operator is of type $(k, 1; a,b)$.
 A type $(1,0;1,0)$ element such as
 $\ad_{\Lambda_2}\overline{\rho}$ could be $\dbar$-exact only if it is identically zero.

By Lemma \ref{image of adjoint} and Formula (\ref{image of adjoint 2}),
$\liet^{*(0,1)}\oplus  \liet^{1,0}\oplus \liec^{1,0}$ is in the kernel of $\ad_{\Lambda_2}$.
If $\ad_{\Lambda_2}{\overline\rho}=0$ in addition, then $\ad_{\Lambda_2}$ vanishes on $\lieg^{*(0,1)}\oplus
\lieg^{1,0}$. It follows that the action of $\ad_{\Lambda_2}$ on $B^{p,q}$
is identically zero for all $p, q\geq 0$.

\begin{lemma}\label{necessary on lambda2}
 If the spectral sequence of the bi-complex of $\Lambda=\Lambda_1+\Lambda_2$ degenerates on the
first page, then $\Lambda_2$ is the center of the Schouten algebra
$\oplus_{p,q}B^{p,q}$.
\end{lemma}

Now we proceed by assuming that $\Lambda_2$ satisfies the necessary condition in the last lemma.
In particular, $\adL=\ad_{\Lambda_1}+\ad_{\Lambda_2}=\ad_{\Lambda_1}$.

Since $\liec^{1,0}$ is one-dimensional, $\Lambda_1=W\wedge T$ for some $T\in \liet^{1,0}$.
$\ad_{\Lambda_1}$ is $\dbar$-exact if there exists
$V$ in $\liet^{1,0}$ such that $\lbra{\Lambda_1}{\overline{\rho}} =\dbar V$.
Under this condition, we compute $\ad_{\Lambda_1}$.

\begin{lemma} Suppose that $\ad_{\Lambda_1}{\overline{\rho}} $ is $\dbar$-exact, then
the map $\ad_{\Lambda_1}$ sends any element in $H^q(\lieg^{p,0})$ to a $\dbar$-exact element.
\end{lemma}
\bproof
Given the dimension constraint on $\liec^{1,0}$, (\ref{image of adjoint 1}) shows that the action of
$\ad_{\Lambda_1}$ is non-trivial only possibly when it acts on components of type
$(k,0; a,1)$. Therefore, if $\Upsilon $ represents a class in $H^q(\lieg^{p,0})$, $\ad_{\Lambda_1}$ could
possibly be non-trivial only when $\Upsilon$ is of type $(p,0;q-1, 1)$.

Given such $\Upsilon$, there exist finitely many
$\overline{\Omega}_j$ in $\liet^{*(0, q-1)}$  and the same number of $\Theta_j$ in
$\liet^{p,0}$ such that
\[
\Upsilon =\overline{\rho}\wedge \sum_{j}(\overline{\Omega}_j\wedge \Theta_j).
\]
Therefore, $\lbra{\Lambda_1}{\Upsilon}$ is equal to
\[
\lbra{\Lambda_1}{\overline{\rho}\wedge \sum_j(\overline{\Omega}_j\wedge \Theta_j)}
=W \wedge \lbra{T}{\overline{\rho}\wedge \sum_j(\overline{\Omega}_j\wedge \Theta_j)}
-T \wedge \lbra{W}{\overline{\rho}\wedge \sum_j(\overline{\Omega}_j\wedge \Theta_j)}.
\]
By Lemma \ref{central vector field}, the last term on the right hand side is identical zero.
By Lemma \ref{adjoint T}, the first term on the right hand side is equal to
\[
W\wedge \lbra{T}{\overline{\rho}}\wedge \sum_j(\overline{\Omega}_j\wedge \Theta_j)
= \lbra{W\wedge T}{\overline{\rho}}\wedge \sum_j(\overline{\Omega}_j\wedge \Theta_j)
\]
If $\ad_{\Lambda_1}\overline{\rho}$ is $\dbar$-exact, then there exists $V$ in $\liet^{1,0}$
such that $\lbra{W\wedge T}{\overline{\rho}}=\lbra{\Lambda_1}{\overline{\rho}}=\dbar V$.
Therefore,
\begin{equation}\label{dbar closed 1}
\lbra{\Lambda_1}{\Upsilon}=
\lbra{W\wedge T}{\overline{\rho}}\wedge \sum_j(\overline{\Omega}_j\wedge \Theta_j)
=\dbar V \wedge (\sum_j(\overline{\Omega}_j\wedge \Theta_j)).
\end{equation}
Given that $\overline\rho$ is $\dbar$-close,
\begin{equation}\label{dbar closed 2}
\dbar\Upsilon=\dbar(\overline{\rho}\wedge  \sum_j(\overline{\Omega}_j\wedge \Theta_j)   )
=-\overline{\rho}\wedge \dbar( \sum_j(\overline{\Omega}_j\wedge \Theta_j) ).
\end{equation}
Since $\dbar(\sum_j(\overline{\Omega}_j\wedge \Theta_j))$ is contained in
$\liet^{*(0,q)}\otimes  \liet^{p-1,0}\otimes \liec^{1,0}$, the exterior product on the right
of the equality (\ref{dbar closed 2}) is equal to zero only if
$\dbar(\sum_j(\overline{\Omega}_j\wedge \Theta_j))=0$. Therefore,
$\Upsilon$ is $\dbar$-closed  if and only if
$\sum_j(\overline{\Omega}_j\wedge \Theta_j)$ is $\dbar$-closed.
It follows that Identity (\ref{dbar closed 1}) is further transformed to
\begin{equation}
\lbra{\Lambda_1}{\Upsilon}=
\dbar (~V \wedge \sum_j(\overline{\Omega}_j\wedge \Theta_j)~ ).
\end{equation}
It shows that the image of any $\dbar$-closed element via $\ad_{\Lambda_1}$ is $\dbar$-exact
when $\ad_{\Lambda_1}\overline\rho$ is $\dbar$-exact.
\eproof

As a consequence of the last two lemmas,
 we conclude that for the nilmanifolds and complex structures as specified, the necessary condition for
the bi-complex of $\Lambda$ to degenerate on the first page is also sufficient.

\begin{theorem} Let $M=G/\Gamma$ be a 2-step nilmanifold with abelian complex structure.
Let $\liec$ be the center of the Lie algebra $\lieg$ of the simply connected Lie group $G$.
Let $\lieg^{1,0}=\liet^{1,0}\oplus \liec^{1,0}$ be the space of invariant $(1,0)$-vectors.
Assume that $\dim\liec^{1,0}=1$.
Suppose that $\Lambda=\Lambda_1+\Lambda_2$ is an invariant holomorphic Poisson bivector with
$\Lambda_1\in  \liet^{1,0}\otimes \liec^{1,0}$ and
$\Lambda_2\in \liet^{2,0}$. Let $\overline\rho$ span $\liec^{*(0,1)}$. The spectral
sequence of the bi-complex of $\Lambda$ degenerates on the first page if and only if
\begin{itemize}
\item $\ad_{\Lambda_1}\overline\rho$ is $\dbar$-exact; and
\item $\Lambda_2$ is in the center of the Schouten algebra $\oplus_{p,q}B^{p,q}$.
\end{itemize}
\end{theorem}

\section{Condition for $\dbar$-exactness}
In this section, we explore deeper in the situation when there exists $V$ such that
\begin{equation}\label{exactness 1}
\lbra{W\wedge T}{\overline{\rho}}=\dbar  V.
\end{equation}
On the left of the equality above,
\begin{equation}\label{exactness 2}
\lbra{W\wedge T}{\overline{\rho}}
=W\wedge  \lbra{ T}{\overline{\rho}}=W\wedge \iota_{T}d{\overline{\rho}}.
\end{equation}
To compute $\dbar V$, we applies the Cartan Formula to evaluate $\dbar V$ on a generic element
$\overline{T}$ in $\liet^{0,1}$ and on $\rho$.
\begin{eqnarray*}
\dbar V(\rho, {\overline{T}}) &=& -V(\lbra{\rho}{\overline{T}})
=V(\lbra{\overline{T}}{\rho})=V(\iota_{\overline{T}}d\rho)\\
&=&d\rho( \overline{T}, V)=-d\rho(V, \overline{T})=-(\iota_Vd\rho)(\overline{T}).
\end{eqnarray*}
Therefore, $\dbar V=-W\wedge\iota_Vd\rho$. Now comparing (\ref{exactness 1}) with
(\ref{exactness 2}), we obtain a rather simple identity below.
\begin{equation}\label{dbar exact equation}
\iota_{T}d{\overline{\rho}}=-\iota_V d\rho.
\end{equation}

Since the complex structure $J$ is abelian, $d\rho$ is a type $(1,1)$-form.
So is $d\overline{\rho}$. We could treat their contractions with elements in
$\liet^{1,0}$  as linear maps.
\begin{equation}
d\rho, \quad d{\overline\rho} \quad : \liet^{1,0}\to \liet^{*(0,1)}.
\end{equation}

Suppose that $d\rho$ has a non-trivial kernel and
$\iota_{T}d{\overline{\rho}}=0$. In such case,
$\ad_{W\wedge T}\overline\rho=0$. However, since $d\overline\omega^j=0$ for all $j$,
$\ad_{W\wedge T}\overline\omega^j=0$. As the complex structure is abelian, then the adjoint action of
$W\wedge T$ on $\lieg^{1,0}$ is identically zero. Therefore, $\ad_{W\wedge T}$ is identically zero
on $B^{p,q}$ for all $p, q\geq 0$. It yields the trivial situation when both $\ad_{\Lambda_1}$ and
$\ad_{\Lambda_2}$ on $B^{p,q}$ are identically zero. In such case, the action of
 $\ad_{\Lambda}$ is trivial and hence as action on $B^{p,q}$, $\dbar_{\Lambda}=\dbar+\ad_{\Lambda}
 =\dbar$. The Poisson cohomology is simply the direct sum of Dolbeault cohomology.
\[
H^n_{\Lambda}(M)=\oplus_{p+q=n} H^q(\lieg^{p,0})=\oplus_{p+q=n} H^q(M, \Theta^p).
\]

On the other hand,  $d\rho(\liet^{1,0})$  is a proper subspace of
$\liet^{*(0,1)}$ if $d\rho$ degenerates. If $T$ is not in the kernel of  $d{\overline\rho}$ but
$\iota_Td{\overline\rho}$ is in the complement of $d\rho(\liet^{1,0})$ then the
equation (\ref{dbar exact equation}) does not necessarily have a solution. In the next section, we will
present an example to demonstrate that such situation does occur.
We summarize our discussion when $d\rho$ degenerates as below.

\begin{theorem}\label{secondary theorem}
Let $M=G/\Gamma$ be a 2-step nilmanifold with abelian complex structure.
Let $\liec$ be the center of the Lie algebra $\lieg$ of the simply connected Lie group $G$.
Let $\lieg^{1,0}=\liet^{1,0}\oplus \liec^{1,0}$ be the space of invariant $(1,0)$-vectors.
Suppose that $\dim_{\mathbb{c}}\liec^{1,0}=1$. Let $W$ span $\liec^{1,0}$ and $\rho$ be its dual.
If $d\overline\rho$ degenerates and $T$ is in its kernel, then
$\Lambda=W\wedge T$ is a holomorphic Poisson
 structure such that the spectral sequence of its associated bi-complex degenerates on
the first page.
\end{theorem}

If $d{\overline\rho}$ is non-degenerate, so would $d\rho$. Therefore, for any $T$ the equation
 (\ref{dbar exact equation}) always has a unique solution.

\begin{theorem}\label{main theorem}
Let $M=G/\Gamma$ be a 2-step nilmanifold with abelian complex structure.
Let $\liec$ be the center of the Lie algebra $\lieg$ of the simply connected Lie group $G$.
Let $\lieg^{1,0}=\liet^{1,0}\oplus \liec^{1,0}$ be the space of invariant $(1,0)$-vectors.
Assume that $\dim\liec^{1,0}=1$.
Suppose that $\Lambda=\Lambda_1+\Lambda_2$ is an invariant holomorphic Poisson bivector with
$\Lambda_1\in  \liet^{1,0}\otimes \liec^{1,0}$ and
$\Lambda_2\in \liet^{2,0}$. Let $\overline\rho$ span $\liec^{*(0,1)}$. The spectral
sequence of the bi-complex of $\Lambda$ degenerates on the second page.
If in addition, the map $d\overline\rho$ is non-degenerate and if
$\Lambda_2$ is in the center of the Schouten algebra $\oplus_{p,q}\lieg^{p,0}\otimes\lieg^{*(0,q)}$,
the spectral sequence degenerates on its first page.
\end{theorem}

By setting $\Lambda_2=0$ in Theorem \ref{main theorem}, we derive
Theorem \ref{corollary theorem} as given in the Introduction.

\section{Examples}

As seen in Theorem \ref{main theorem}, invariant holomorphic Poisson
structure $\Lambda$ has two components $\Lambda_1$ and $\Lambda_2$. If we requires the spectral sequence to
degenerate on the first page,
the action of $\ad_{\Lambda_2}$ is identically zero. Therefore, in this section
we focus on the case when $\Lambda_2=0$, and hence
$\Lambda=\Lambda_1=W\wedge T$ where $W \in\liec^{1,0}$ and $T\in \liet^{1,0}$.
 As the complex structure is abelian, it is obvious that
 $\lbra{W\wedge T}{W\wedge T}=0$. In addition, since $\dbar W=0$, $\dbar\Lambda=-W\wedge \dbar T$.
 By (\ref{adj-T-on-rho-bar}),  $\dbar T$ is in $\liec^{1,0}\otimes \liet^{*(0,1)}$. When
 $\dim\liec^{1,0}=1$, $W\wedge \dbar T=0$. Therefore, for any $T\in \liet^{1,0}$,
 the bivector $\Lambda=W\wedge T$ is an invariant holomorphic Poisson structure.
 We now focus on this type of structures, and present concrete examples to illustrate
 Theorem \ref{secondary theorem} as well as Theorem \ref{main theorem}.

 \

 When $\dim\liec^{1,0}=1,$  the dual structure equations in
 (\ref{dual 1}) and (\ref{dual 2})
 are simplified as below.
\[
d\rho=\sum_{i,j}E_{ji}\omega^i\wedge\oomega^j,
\quad \mbox{ and } \quad
d{\overline{\rho}}
=-\sum_{i,j}{\overline{E}_{ij}}\omega^i\wedge\oomega^j.
\]

\noindent{\bf Example 1.} Consider a one-dimensional central extension of the Heisenberg
algebra $\lieh_{2n+1}$ of real dimension $2n+1$. Let $\{X_j, Y_j, Z, A: 1\leq j\leq n\}$
be basis with structure equation
\begin{equation}
\lbra{X_j}{Y_j}=-\lbra{Y_j}{X_j}=Z, \quad \mbox{ for all } \quad 1\leq j\leq n.
\end{equation}
The real center $\liec$ is spanned by $Z$ and $A$. Define an almost complex structure by
\[
JX_j=Y_j, \quad JY_j=-X_j, \quad JZ=A, \quad JA=-Z.
\]
It is an abelian complex structure. Let $W=\frac12(Z-iA)$ and $T_j=\frac12(X_j-iY_j)$, then the complex structure
equation becomes
\[
\lbra{{\overline{T}_j}}{T_j}=-\frac{i}{2}(W+\overline{W}).
\]
Therefore, $E_{ii}=-\frac{i}{2}=-\overline{E}_{ii}$. Hence $d\rho$ is non-degenerate and serves
as an example for Theorem \ref{main theorem}.

\

\noindent{\bf Example 2.} On the direct sum of two Heisenberg algebras $\lieh_{2m+1}\oplus \lieh_{2n+1}$, choose
a basis $\{X_j, Y_j, Z, A_k, B_k, C; 1\leq j\leq m, 1\leq k\leq n\}$
so that the non-zero structure equations are give by
\begin{equation}
\lbra{X_j}{Y_j}=-\lbra{Y_j}{X_j}=Z, \quad \lbra{A_k}{B_k}=-\lbra{B_k}{A_k}=C.
\end{equation}
Define an almost complex structure $J$ by
\[
JX_j=Y_j, \quad JA_k=B_k, \quad JZ=C.
\]
We get an abelian complex structure. Let
\[
W=\frac12(Z-iC), \quad S_j=\frac12(X_j-iY_j), \quad T_k=\frac12(A_k-iB_k).
\]
$\liec^{1,0}$ is spanned by $W$.
It follows that the non-zero complex structure equations are given as below.
\[
\lbra{\overline{S}_j}{S_j}=-\frac{i}2(W+{\overline{W}}), \quad
\lbra{\overline{T}_k}{T_k}=\frac12(W-{\overline{W}}).
\]
Then the structural constants
\[
E_{jj}=-\frac{i}2=-{\overline{E}_{jj}},
\quad \mbox{ and } \quad
E_{kk}=\frac12 ={\overline{E}_{kk}} ,
\]
for all $1\leq j\leq m$ and $1\leq k\leq n$.
It is now obvious that $d\rho$ is non-degenerate and serves as example for
Theorem \ref{main theorem}.

\

\noindent{\bf Example 3.} Consider a real vector space $W_{4n+6}$ spanned by \[
\{X_{4k+1}, X_{4k+2}, X_{4k+3},
X_{4k+4}, Z_1, Z_2; 0\leq k\leq n\}.\] Define a Lie bracket by
\begin{eqnarray*}
\lbra{X_{4k+1}}{X_{4k+3}}=-\frac12Z_1,
& \quad &
\lbra{X_{4k+1}}{X_{4k+4}}=-\frac12Z_2, \\
\lbra{X_{4k+2}}{X_{4k+3}}=-\frac12Z_2,
& \quad &
\lbra{X_{4k+2}}{X_{4k+3}}=\frac12Z_1.
\end{eqnarray*}
One could define an abelian complex structure $J$ by
\[
JX_{4k+1}=X_{4k+2}, \quad JX_{4k+3}=-X_{4k+4}, \quad JZ_1=-Z_2.
\]
For $0\leq k\leq n$, define
\[
W=\frac12(Z_1+iZ_2), \quad
T_{2k+1}=\frac12 (X_{4k+1}-iX_{4k+2}),
\quad
T_{2k+2}=\frac12 (X_{4k+3}+iX_{4k+4}).
\]
It is now a straightforward computation to show that the non-zero structure equations are
\begin{equation}
\lbra{{\overline{T}}_{2k+1}}{T_{2k+2}}=-\frac12 W, \quad
\lbra{{\overline{T}}_{2k+2}}{T_{2k+1}}-=\frac12{\overline{W}}.
\end{equation}
Except when
\[
E_{2k+1,2k+2}=-\frac12, \quad \mbox{ for all } \quad 0\leq k\leq n,
\]
all other structure constants are equal to zero. In particular,
\[
d\rho=-\frac12\sum_{k=0}^n\omega^{2k+2}\wedge\oomega^{2k+1}, \quad
\mbox{ and }
\quad
d\overline{\rho}=\frac12 \sum_{k=0}^n\omega^{2k+1}\wedge\oomega^{2k+2}.
\]
Treating $d\rho$ and $d\overline{\rho}$ as maps from $\liet^{1,0}$ to
$\liet^{0,1}$, their image spaces are transversal.
Therefore, given $T\in \liet^{1,0}$ such that there exists $V\in \liet^{1,0}$ with
$\iota_Td\overline{\rho}=-\iota_Vd{\rho}$ only if $\iota_Td\overline{\rho}=0$. It is possible only when
$T$ is the complex linear span of
$\{ T_2, \dots, T_{2k+2}, \dots, T_{2n+2} \}$, which is the kernel of $d\overline\rho$.

In this case, it simply means that $\lbra{T}{\overline\rho}=0$, and hence $\adL=0$. This example illustrates the
situation in Theorem \ref{secondary theorem}.

\

\noindent{Example 4.} Consider a real vector space $P_{4n+2}$ spanned by
\[
\{X_{4k+1}, X_{4k+2}, X_{4k+3},
X_{4k+4}, Z_1, Z_2; 0\leq k\leq n\}.\] Define a Lie bracket by
\begin{eqnarray*}
\lbra{X_{4k+1}}{X_{4k+2}}=-\frac12Z_1,
\quad
\lbra{X_{4k+1}}{X_{4k+4}}=-\frac12Z_2,
\quad
\lbra{X_{4k+2}}{X_{4k+3}}=-\frac12Z_2.
\end{eqnarray*}
Define an abelian complex structure $J$ on  $P_{4n+2}$ by
\[
J{X_{4k+1}}={X_{4k+2}},
\quad
J{X_{4k+3}}=-{X_{4k+4}},
\quad
J{Z_1}=-Z_2.
\]
Define
\[
W=\frac12(Z_1+iZ_2), \quad
T_{2k+1}=\frac12 (X_{4k+1}-iX_{4k+2}),
\quad
T_{2k+2}=\frac12 (X_{4k+3}+iX_{4k+4}).
\]
The non-zero structure equations in terms of complex vectors are
\[
\lbra{{\overline{T}}_{2k+1}}{T_{2k+1}}=\frac{i}4(W+\overline{W});
\quad
\lbra{{\overline{T}}_{2k+1}}{T_{2k+2}}=-\frac14(W-\overline{W}).
\]
Hence
\[
E_{2k+1, 2k+1}=\frac{i}4,
\quad
E_{2k+1, 2k+2}=-\frac14,
\quad
E_{2k+2, 2k+1}=-\frac14.
\]
It follows that $d\rho$ is non-degenerate, and hence (\ref{dbar exact equation}) is solvable.
We have an example for Theorem \ref{main theorem} again.

\

\noindent{\bf Acknowledgment.}
Y.~S. Poon thanks the Yau Mathematical Sciences Center at
Tsinghua University for hospitality where he finished this paper in summer 2016.

\end{document}